\documentclass{amsart} 
\usepackage{amssymb,upref,enumerate}
\newtheorem{theorem}{Theorem}[section]
\newtheorem{lemma}[theorem]{Lemma}
\newtheorem{proposition}[theorem]{Proposition}
\newtheorem{corollary}[theorem]{Corollary}
\theoremstyle{definition}
\newtheorem{definition}[theorem]{Definition}
\newtheorem{example}[theorem]{Example}
\theoremstyle{remark}
\newtheorem{remark}[theorem]{Remark}

\newtheorem{question}[theorem]{Question}
\numberwithin{equation}{section}
\DeclareMathOperator{\diam}{diam}
\DeclareMathOperator{\dist}{dist}

\DeclareMathOperator*{\essinf}{ess\,inf\,}
\newcommand{\abs}[1]{\lvert #1 \rvert}
\newcommand{\Abs}[1]{\left\lvert #1 \right\rvert}
\newcommand{\norm}[1]{\lVert #1 \rVert}
\newcommand{\inn}[1]{\langle #1 \rangle}

\newcommand{\LL}{\mathcal{L}}
\newcommand{\R}{\mathbb{R}}
\newcommand{\Ri}{\mathcal{I}}
\newcommand{\Z}{\mathbb{Z}}

\begin{document}

\title{Doubling measures, monotonicity, and quasiconformality}

\author{Leonid V. Kovalev}
\address{Department of Mathematics, Mailstop 3368, Texas A\&M University,
College Station, TX 77843-3368 \\ USA}
\email{lkovalev@math.tamu.edu}
\thanks{L.V.K. was supported by an NSF Young Investigator award
 under grant DMS 0601926.}
\author{Diego Maldonado}
\address{Department of Mathematics \\ University of Maryland \\
College Park, MD 20742 \\ USA}
\email{maldona@math.umd.edu}
\author{Jang-Mei Wu}
\address{Department of Mathematics \\ University of Illinois \\
1409 West Green Street \\ Urbana, IL 61801 \\ USA}
\email{wu@math.uiuc.edu}
\thanks{J.-M.W. was supported by the NSF grant DMS 0400810.}

\begin{abstract} We construct quasiconformal mappings in Euclidean spaces by
integration of a discontinuous kernel against doubling measures with suitable decay.
The differentials of mappings that arise in this way satisfy an isotropic form of the doubling condition.
We prove that this isotropic doubling condition is satisfied by the distance functions
of certain fractal sets. Finally, we construct an isotropic doubling measure that is
not absolutely continuous with respect to the Lebesgue measure.
\end{abstract}

\maketitle

\section{Introduction}

Given a nonatomic positive Radon measure $\mu$ on $\R$, define an increasing function
$f_{\mu}(x)=\int_0^x d\mu(z)$, so that $f_{\mu}'=\mu$ in the sense of distributions.
It is well-known that $\mu$ is doubling if and only if $f_{\mu}$ is quasisymmetric.
We extend this relation between doubling measures and quasisymmetric mappings
to higher dimensions. Observe that
\[
f_{\mu}(x)=\frac{1}{2}\int_{\R}\left\{\frac{x-z}{\abs{x-z}}+\frac{z}{\abs{z}}\right\}\,d\mu(z).
\]
For a nonatomic Radon measure $\mu$ on $\R^n$, $n\ge 1$, we define $f_{\mu} : \R^n\to\R^n$ by
\begin{equation}\label{db0a}
f_{\mu}(x)=\frac{1}{2}\int_{\R^n}\left\{\frac{x-z}{\abs{x-z}}+\frac{z}{\abs{z}}\right\}\,d\mu(z).
\end{equation}
where $\abs{\cdot}$ is now interpreted as the Euclidean norm.
When $n\ge 2$ the integrand in~\eqref{db0a} is not compactly supported with
respect to $z$. The integral~\eqref{db0a} converges provided that $\mu$ satisfies the decay condition
\begin{equation}\label{db4}
\int_{\abs{z}>1}\abs{z}^{-1}\,d\mu(z)<\infty.
\end{equation}
For $0<\gamma<n$, let
\[\Ri_{\gamma}\mu(x)=\int_{\R^n}\abs{x-z}^{\gamma-n}\,d\mu(z)\]
be the Riesz potential of $\mu$ of
order $\gamma$. Condition~\eqref{db4} is equivalent to $\Ri_{n-1}\mu$ being finite almost everywhere.
Here and in the sequel the words ``almost everywhere'' or ``a.e.'' refer to the Lebesgue measure.
A positive Radon measure $\mu$ on $\R^n$ is called doubling if there is a constant $C\ge 1$ such that
$\mu(2B)\le C\mu(B)$ for every ball $B\subset\R^n$. Here $2B$ stands for the ball that has the same center
as $B$ and twice its radius.

\begin{theorem}\label{dbthm} Let $\mu$ be a doubling measure on $\R^n$ that satisfies the
decay condition~\eqref{db4}.
The mapping $f_{\mu}$ defined by~\eqref{db0a} is $\eta$-quasisymmetric
with $\eta$ depending only on the doubling constant of $\mu$.
Furthermore, $f_{\mu}$ is $\delta$-monotone and for a.e. $x\in\R^n$
\begin{equation}\label{Jacdo1}
\frac{1}{12C^3}\Ri_{n-1}\mu(x)\le \norm{Df_{\mu}(x)}\le \pi\Ri_{n-1}\mu(x),
\end{equation}
where $C$ is the doubling constant of $\mu$.
\end{theorem}
Theorem~\ref{dbthm} expands the class of weights that are known to be comparable to
Jacobians of quasiconformal mappings. The quasiconformal Jacobian problem posed by David and
Semmes in~\cite{DS90} asks for a characterization of all such weights (see also~\cite{Bi,BHS04,BHS06,La02,Se96}).
The authors of~\cite{BHS04} point out that such a characterization would give a good idea
of which metric spaces are bi-Lipschitz equivalent to $\R^n$.

A mapping $f\colon\R^n\to\R^n$ is called monotone if
$\inn{F(x)-F(y),x-y}\ge 0$ for all $x,y\in\R^n$, where $\inn{\cdot,\cdot}$ is the inner product~\cite{AA99,Ze90}.
In other words, $F$ is monotone if the angle formed
by the vectors $F(x)-F(y)$ and $x-y$ is at most $\pi/2$. A stronger version
of this condition, introduced by Sobolevskii in~\cite{So57},
requires the angle to be bounded by a constant less than $\pi/2$.

\begin{definition}\label{defdm} Let $\delta\in (0,1]$.
A mapping $F$ from a convex domain $\Omega\subset \R^n$ into $\R^n$ is
called $\delta$-monotone if for all $x,y\in\Omega$
\begin{equation}\label{deflr}
\inn{F(x)-F(y),x-y}\ge\delta\abs{F(x)-F(y)}\abs{x-y}.
\end{equation}
\end{definition}
Let $\Omega$ be a domain in $\R^n$. An injective mapping $f : \Omega\to\R^n$ is called
$\eta$-quasi\-symmetric if there is
a homeomorphism $\eta : [0,\infty)\to [0,\infty)$ such that
\begin{equation}\label{def3a}
\frac{\abs{f(x)-f(z)}}{\abs{f(y)-f(z)}}\le
\eta\left(\frac{\abs{x-z}}{\abs{y-z}}\right),\quad z\in \Omega, \
x,y\in\Omega\setminus\{z\}.
\end{equation}
In the case $n\ge 2$, $f$ is called $K$-quasiconformal
if $f\in W_{\rm loc}^{1,n}(\Omega;\R^n)$ and the operator norm of the derivative $Df$ satisfies
$\norm{Df(x)}^n\le K \det Df(x)$, a.e. in $\Omega$ for some $K\ge 1$.

Every sense-preserving quasisymmetric mapping $f\colon \R^n\to\R^n$, $n\ge 2$, is quasiconformal
and vice versa (\cite[Ch.~10]{He01}, \cite[p.~98]{TV80}). The following result relates $\delta$-monotone
and quasisymmetric mappings.

\begin{theorem}\label{London}~\cite[Theorem~4]{Ko} Suppose that $n\ge 2$ and
$f\colon \Omega\to\R^n$ is a nonconstant $\delta$-monotone mapping. If $B$ is a closed ball
such that $2B\subset \Omega$, then $f$ is $\eta$-quasisymmetric on $B$ with $\eta$ depending
only on $\delta$.
\end{theorem}

Quasiconformality is known to be related to the doubling condition in several ways~\cite{BHR01,Ro01,St92}.
For example, if $f\colon\R^n\to\R^n$ is quasiconformal, then $\norm{Df}^n$ is
a doubling weight, and moreover an $A_{\infty}$ weight~\cite{Ge73}. Consequently, $\norm{Df}$ is doubling as well.
Since $\delta$-monotonicity is a stronger property than quasiconformality, one can expect that the differential
of a $\delta$-monotone mappings exhibits a stronger doubling behavior.
Our Theorem~\ref{dmidc} confirms this. Before stating it, we observe that a Radon measure $\mu$ is doubling
if and only if there exists a constant $A$ such that
\[
A^{-1}\le \frac{\mu(Q_1)}{\mu(Q_2)}\le A
\]
for any congruent cubes $Q_1$ and $Q_2$ with nonempty intersection.
Recall that two subsets of $\R^n$ are called congruent if there is an isometry of $\R^n$
that maps one of them onto the other.

\begin{definition}\label{idmdef} A Radon measure $\mu$ on $\R^n$ is \emph{isotropic doubling}
if there is a constant $A\ge 1$ such that
\begin{equation}\label{idmdef1}
A^{-1}\le \frac{\mu(R_1)}{\mu(R_2)}\le A
\end{equation}
whenever $R_1$ and $R_2$ are congruent rectangular boxes with nonempty intersection.
\end{definition}

\begin{theorem}\label{dmidc} Let $f : \R^n\to\R^n$ be a nonconstant $\delta$-monotone mapping, $n\ge2$.
Then the weight $\norm{Df}$ is isotropic doubling.
\end{theorem}

The requirement~\eqref{idmdef1} for all boxes regardless of
their orientation and aspect ratio imposes a very strong condition on the measure when $n\ge 2$.
In particular, in any cube the projection of an isotropic doubling measure
to any $(n-1)$-dimensional face of the cube is comparable to the
$(n-1)$-dimensional Lebesgue measure $\LL^{n-1}$ (Lemma~\ref{idc1}).

\begin{theorem}\label{sing} For every $n\ge 2$ there exists an isotropic doubling
measure $\mu$ on $\R^n$ that is purely singular with respect to $\LL^n$.
\end{theorem}

\begin{theorem}\label{nobl} For every $n\ge 2$ there exists an isotropic doubling measure $\mu$ on $\R^n$ and a
bi-Lipschitz mapping $f\colon\R^n\to\R^n$ such that the pushforward measure $f_{\#}\mu$ is not isotropic doubling.
\end{theorem}

Theorems~\ref{dbthm}, \ref{dmidc}, \ref{sing}, and \ref{nobl} are proved in sections~\ref{doub}, \ref{idm},
\ref{singular} and \ref{BL}, respectively.
We also prove that the distance functions of certain fractal subsets of $\R^n$ give
rise to isotropic doubling weights (Proposition~\ref{ULNC}).
By virtue of this result, the sets constructed by
Semmes~\cite{Se96} and Laakso~\cite{La02} provide examples of isotropic doubling weights
that are not comparable to $\norm{Df}$ for any $\delta$-monotone, or even quasiconformal,
mapping $f\colon\R^n\to\R^n$.

\section{Doubling measures and monotone mappings}\label{doub}

The balls, cubes, and rectangular boxes considered in this paper are assumed closed.
A nonnegative locally integrable function on $\R^n$ is called a \emph{weight}.
A weight is doubling if the measure $w(x)d\LL^n(x)$ is doubling, where $\LL^n$ is the $n$-dimensional
Lebesgue measure. Throughout the paper we only consider positive nonzero measures.

In this section we study the mapping $f_{\mu}$ defined by~\eqref{db0a}. In particular,
we prove Theorem~\ref{dbthm} and its more general version, Theorem~\ref{dbthm2}.
Given three distinct points $x,y,z\in\R^n$,
let $l(x,y,z)=\abs{x-y}+\abs{x-z}+\abs{y-z}$ denote the perimeter of the triangle $xyz$. Also let
$\widehat{zxy}$ be the angle between the vectors $y-x$ and $z-x$, and define
\[
\tau(x,y,z):=\pi-\max\{\widehat{zxy},\widehat{zyx}\}.
\]

\begin{theorem}\label{dbthm2} Let $\Omega\subset\R^n$ be a convex domain, $n\ge2$.
Let $\mu$ be a nonatomic Radon measure on $\R^n$ that satisfies~\eqref{db4}.
Suppose that there exists $\kappa>0$ such that
\begin{equation}\label{dbthm3}
\liminf_{y\to x} \left(\int_{\R^n}\frac{\tau^2(x,y,z)\,d\mu(z)}{l(x,y,z)}\right)
\bigg/ \left(\int_{\R^n}\frac{\tau(x,y,z)\,d\mu(z)}{l(x,y,z)}\right)\ge \kappa
\end{equation}
for all $x\in\Omega$.
Then $f_{\mu}$ is $\delta$-monotone in $\Omega$ with $\delta=\kappa/(2\pi^2)$.
\end{theorem}

The proof is preceded by an elementary lemma.

\begin{lemma}\label{dblem} Fix $z\in\R^n$ and define $g_z(x)=(x-z)/\abs{x-z}$ for $x\in\R^n\setminus\{z\}$.
For any distinct points $x,y\in\R^n\setminus\{z\}$ we have
\begin{equation}\label{dblem1m}
\frac{2\abs{x-y}}{\pi l(x,y,z)}\tau(x,y,z)\le\abs{g_z(x)-g_z(y)}\le \frac{4\abs{x-y}}{l(x,y,z)}\tau(x,y,z)
\end{equation}
and
\begin{equation}\label{dblem2m}
\frac{2\abs{x-y}^2}{\pi^2 l(x,y,z)}\tau(x,y,z)^2\le \inn{g_z(x)-g_z(y),x-y}\le
\frac{4\abs{x-y}^2}{l(x,y,z)}\tau(x,y,z)^2.
\end{equation}
\end{lemma}
\begin{proof} First we prove~\eqref{dblem1m}. By the sine theorem
\begin{equation}\label{dblem1a}
\frac{\sin\widehat{xzy}}{\abs{x-y}}=\frac{\sin\widehat{zyx}+\sin\widehat{zxy}}{\abs{x-z}+\abs{y-z}}.
\end{equation}
Express the sum of sines as a product:
\begin{equation}\label{dblem1b}
\sin\widehat{zyx}+\sin\widehat{zxy}=2\cos\frac{\widehat{xzy}}{2}\cos\frac{\widehat{zyx}-\widehat{zxy}}{2}.
\end{equation}
Combining~\eqref{dblem1a}, ~\eqref{dblem1b}, and the identity
$\sin\widehat{xzy}=2\sin\frac{\widehat{xzy}}{2}\cos\frac{\widehat{xzy}}{2}$, we obtain
\begin{equation}\label{dblem1c}
\frac{\sin\frac{\widehat{xzy}}{2}}{\abs{x-y}}=\frac{\cos\frac{\widehat{zyx}-\widehat{zxy}}{2}}{\abs{x-z}+\abs{y-z}}.
\end{equation}
The definition of $g_z$ implies
\[
\abs{g_z(x)-g_z(y)}=2\sin\frac{\widehat{xzy}}{2},
\]
which together with~\eqref{dblem1c} yield
\begin{equation}\label{dblem1}
\abs{g_z(x)-g_z(y)}= \frac{2\abs{x-y}}{\abs{x-z}+\abs{y-z}}\cos\frac{\widehat{zyx}-\widehat{zxy}}{2}.
\end{equation}
By the triangle inequality
\begin{equation}\label{lest}
\frac12 l(x,y,z)\le \abs{x-z}+\abs{y-z}\le l(x,y,z).
\end{equation}
Therefore~\eqref{dblem1m} will follow from~\eqref{dblem1} once we prove that
\begin{equation}\label{dblem1d}
\frac{1}{\pi}\tau(x,y,z)\le \cos\frac{\widehat{zyx}-\widehat{zxy}}{2}\le \tau(x,y,z).
\end{equation}
Let $\alpha=\max\{\widehat{zxy},\widehat{zyx}\}$. We have
\[\cos\frac{\widehat{zyx}-\widehat{zxy}}{2}\ge \cos\frac{\alpha}{2}=
\sin\frac{\pi-\alpha}{2}\ge \frac{2}{\pi}\frac{\pi-\alpha}{2}=\frac{1}{\pi}\tau(x,y,z),
\]
which proves the first inequality in~\eqref{dblem1d}. The second inequality is
trivial when $\alpha\le \pi/2$, because then $\tau(x,y,z)=\pi-\alpha>1$.
If $\alpha>\pi/2$, then
\[
\cos\frac{\widehat{zyx}-\widehat{zxy}}{2}\le \cos\frac{\alpha-(\pi-\alpha)}{2}
=\sin(\pi-\alpha)\le \pi-\alpha=\tau(x,y,z).
\]
This proves~\eqref{dblem1d} and ~\eqref{dblem1m}.

The proof of~\eqref{dblem2m} involves the identity
$\cos\alpha+\cos\beta=2\cos\frac{\alpha+\beta}{2}\cos\frac{\alpha-\beta}{2}$ as follows:
\begin{equation*}\begin{split}
\inn{g_z(x)-g_z(y),x-y}&=\frac{\inn{z-x,y-x}}{\abs{x-z}}+\frac{\inn{z-y,x-y}}{\abs{y-z}}\\
&=(\cos\widehat{zxy}+\cos\widehat{xyz})\abs{x-y}\\
&=2\sin\frac{\widehat{xzy}}{2}\cos\frac{\widehat{zxy}-\widehat{xyz}}{2}\abs{x-y}\\
&=\abs{g_z(x)-g_z(y)}\abs{x-y}\cos\frac{\widehat{zxy}-\widehat{xyz}}{2}.
\end{split}
\end{equation*}
Applying~\eqref{dblem1} we obtain
\begin{equation}\label{dblem2}
\inn{g_z(x)-g_z(y),x-y}=\frac{2\abs{x-y}^2}{\abs{x-z}+\abs{y-z}}\cos^2\frac{\widehat{zyx}-\widehat{zxy}}{2}.
\end{equation}
This together with~\eqref{lest} and~\eqref{dblem1d} imply~\eqref{dblem2m}.
\end{proof}

\begin{proof}[Proof of Theorem~\ref{dbthm2} ] Since
\[f_{\mu}(x)-f_{\mu}(y)=\frac12\int_{\R^n}(g_z(x)-g_z(y))\,d\mu(z),\]
Lemma~\ref{dblem} implies the following inequalities.
\begin{equation}\label{do6}
\abs{f_{\mu}(x)-f_{\mu}(y)}\le 2\abs{x-y}\int_{\R^n}\frac{\tau(x,y,z)}{l(x,y,z)}d\mu(z);
\end{equation}
\begin{equation}\label{do3}
\begin{split}
\frac{1}{\pi^2}\abs{x-y}^2\int_{\R^n} \frac{\tau(x,y,z)^2}{l(x,y,z)}\,d\mu(z)
&\le \inn{f_{\mu}(x)-f_{\mu}(y),x-y} \\ &\le
2\abs{x-y}^2\int_{\R^n} \frac{\tau(x,y,z)^2}{l(x,y,z)}\,d\mu(z).
\end{split}
\end{equation}
Choose a number $\kappa'$ so that $0<\kappa'<\kappa$.
Combining~\eqref{do6}, \eqref{do3}, and~\eqref{dbthm3},
we conclude that for every $x\in\Omega$ there exists $\varepsilon(x)>0$ such that
\begin{equation}\label{do8}
\inn{f_{\mu}(x)-f_{\mu}(y),x-y}\ge \frac{\kappa'}{2\pi^2}\abs{f_{\mu}(x)-f_{\mu}(y)}\abs{x-y}
\end{equation}
whenever $\abs{x-y}\le\varepsilon(x)$.
Next, let $x$ and $y$ be any distinct points in $\Omega$. The line segment $[x,y]$
is covered by open balls $B(z,\varepsilon(z)/2)$, $z\in[x,y]$. Choose a finite subcover
with centers $z_j=y+t_j(x-y)$, $0=t_0<t_1<\dots<t_N=1$. Then
\[\begin{split}
\inn{f_{\mu}(x)-f_{\mu}(y),x-y}&=\sum_{j=1}^N \inn{f_{\mu}(z_j)-f_{\mu}(z_{j-1}),x-y}\\
&\ge \frac{\kappa'}{2\pi^2} \sum_{j=1}^N \abs{f_{\mu}(z_j)-f_{\mu}(z_{j-1})}\abs{x-y}\\
&\ge \frac{\kappa'}{2\pi^2} \abs{f_{\mu}(x)-f_{\mu}(y)}\abs{x-y}.
\end{split}\]
Letting $\kappa'\to\kappa$ completes the proof.
\end{proof}

\begin{proof}[Proof of Theorem~\ref{dbthm}] Let $C$ be such that $\mu(2B)\le C\mu(B)$ for any ball $B\subset\R^n$.
Given two distinct points $x,y\in\R^n$ and $r>0$, choose a point $w\in\R^n$ so that
$x-w$ is orthogonal to $x-y$ and $\abs{x-w}=r/2$. It is easy to see that
$\tau(x,y,z)\ge \pi/3$ for all $z\in B(w,r/4)$.
Since $B(x,r)\subset B(w,3r/2)\subset 8B(w,r/4)$, it follows that
$\mu(B(x,r))\le C^3\mu(B(w,r/4))$. Using this together with the inequality
\[
\abs{x-z}+\abs{x-y}\le l(x,y,z)\le 2(\abs{x-z}+\abs{x-y}),
\]
we obtain
\begin{equation}\label{dtr1}
\int_{\R^n}\frac{d\mu(z)}{l(x,y,z)}\le \int_{\R^n}\frac{d\mu(z)}{\abs{x-z}+\abs{x-y}} =
\int_0^{\infty}\frac{\mu(B(x,r))\,dr}{(r+\abs{x-y})^2}
\end{equation}
and
\begin{equation}\label{dtr2}
\int_{\R^n}\frac{\tau(x,y,z)^2}{l(x,y,z)}d\mu(z)\ge
\frac{\pi^2}{6C^3}\int_0^{\infty}\frac{\mu(B(x,r))\,dr}{(r+\abs{x-y})^2}.
\end{equation}
Therefore, $\mu$ satisfies a stronger condition than~\eqref{dbthm3}, namely
\begin{equation}\label{dbthms}
\int_{\R^n}\frac{\tau^2(x,y,z)}{l(x,y,z)}d\mu(z)
\ge c\int_{\R^n}\frac{d\mu(z)}{l(x,y,z)}
\end{equation}
for any distinct points $x,y\in\R^n$. Here $c=\pi^2/(6C^3)$.
By Theorem~\ref{dbthm2} the mapping $f_{\mu}$ is $\delta$-monotone,
and therefore $\eta$-quasisymmetric. Here $\eta$ depends only on $\delta$, which in turn depends only on $C$.

It remains to prove~\eqref{Jacdo1}. Since $l(x,y,z)\ge 2\abs{x-z}$, inequality~\eqref{do6} implies
\begin{equation}\label{Jacup}
\limsup_{y\to x}\frac{\abs{f_{\mu}(x)-f_{\mu}(y)}}{\abs{x-y}}\le \pi\Ri_{n-1}\mu(x),\quad x\in\Omega.
\end{equation}
Inequalities~\eqref{do3} and~\eqref{dbthms} yield
\[
\liminf_{y\to x}\frac{\abs{f_{\mu}(x)-f_{\mu}(y)}}{\abs{x-y}}\ge
\frac{1}{6C^3}\liminf_{y\to x}\int_{\R^n}\frac{d\mu(z)}{l(x,y,z)}.
\]
For any $\varepsilon>0$ the limit
\[\lim_{y\to x}\frac{1}{l(x,y,z)}=\frac{1}{2\abs{x-z}}\]
is uniform with respect to $z\in\R^n\setminus B(x,\varepsilon)$. Therefore,
\[\liminf_{y\to x}\int_{\R^n}\frac{d\mu(z)}{l(x,y,z)} \ge
\int_{\R^n\setminus B(x,\varepsilon)}\frac{d\mu(z)}{2|x-z|}\]
Letting $\varepsilon\to 0$ yields
\begin{equation}\label{Jacdown}
\liminf_{y\to x}\frac{\abs{f_{\mu}(x)-f_{\mu}(y)}}{\abs{x-y}}\ge \frac{1}{12C^3}\Ri_{n-1}\mu(x), \quad x\in\Omega,
\end{equation}
because $\mu$ has no atoms. Combining the estimates~\eqref{Jacup} and~\eqref{Jacdown},
we conclude that~\eqref{Jacdo1} holds whenever $f_{\mu}$ is differentiable at $x$.
\end{proof}

\begin{remark}\label{close} Condition~\eqref{dbthm3} is close to being best possible for the $\delta$-monotonicity
of $f_{\mu}$ in Theorem~\ref{dbthm2} in view of~\eqref{dblem1m}, \eqref{dblem2m}, and \eqref{do3}.
The proof of Theorem~\ref{dbthm} shows that condition~\eqref{dbthms} is sufficient for
the validity of~\eqref{Jacdo1}.
\end{remark}

Since the assumption~\eqref{dbthm3} in Theorem~\ref{dbthm2} is hard to verify directly, we
state a simpler condition that implies it. Given a line $L$ in $\R^n$ and a point $x\in L$, let
$S_{\alpha}(x,L)$ be the double-sided cone of opening angle $\alpha\in(0,\pi/2)$ with vertex
$x$ and axis $L$. Formally, $z\in S_{\alpha}(x,L)$ if and only if the acute angle between the vector
$z-x$ and the line $L$ is at most $\alpha$. We also introduce a notation for spherical shells
$A(x,r,R)=\{z\colon r< \abs{z-x}\le R\}$.

\begin{proposition}\label{geom} Let $\mu$ be a nonatomic Radon measure on $\R^n$. Suppose that
there exist constants $\alpha\in (0,\pi/2)$, $C>0$, and $M\ge 1$
such that for any line $L\subset\R^n$, any $x\in L$ and any $r>0$ we have
\begin{equation}\label{geom1}
\mu(A(x,r,2r))\le C\mu (A(x,M^{-1}r,2Mr)\setminus S_{\alpha}(x,L)).
\end{equation}
Then the measure $\mu$ satisfies~\eqref{dbthms} and consequently~\eqref{dbthm3}.
\end{proposition}
\begin{proof} Without loss of generality we may assume that $M=2^m$ for some integer $m$.

Given $y\in\R^n\setminus\{x\}$, let $L$ be the line through $x$ and $y$.
Inequality~\eqref{geom1} implies
\begin{equation}\label{orgcone}
\int_{\R^n}\frac{d\mu(z)}{l(x,y,z)}\le
C_1\int_{\R^n\setminus S_{\alpha}(x,L)}\frac{d\mu(z)}{l(x,y,z)}
\end{equation}
where $C_1$ depends only on $C$ and $M$. Indeed,
\[\begin{split}
\int_{\R^n}\frac{d\mu(z)}{l(x,y,z)}
&=\sum_{k\in\Z}\int_{A(x,2^k,2^{k+1})}\frac{d\mu(z)}{l(x,y,z)}
\le \sum_{k\in\Z}\frac{\mu(A(x,2^k,2^{k+1}))}{2^k+\abs{x-y}}\\
&\le C\sum_{k\in\Z}\frac{\mu(A(x,M^{-1}2^k,M2^{k+1})\setminus S_{\alpha}(x,L))}{2^k+\abs{x-y}}\\
&\le C_1\int_{\R^n\setminus S_{\alpha}(x,L)}\frac{d\mu(z)}{l(x,y,z)},
\end{split}\]
which proves~\eqref{orgcone}.
Since $\tau(x,y,z)\ge\alpha$ for all $z\notin S_{\alpha}(x,L)$, estimate~\eqref{orgcone}
implies~\eqref{dbthms}.
\end{proof}

Any doubling measure satisfies~\eqref{geom1}, because the set $A(x,r,2r)\setminus S_{\alpha}(x,L)$
contains a ball of radius comparable to $r$. The following result shows that~\eqref{geom1}
is satisfied by some non-doubling measures as well.

\begin{proposition}\label{geomI} Let $\nu$ be a doubling measure on $\R^n$, $n\ge 2$.
Define a measure $\mu$ by setting $\mu(E)=\nu(E\setminus B(0,1))$.
Then $\mu$ satisfies the assumptions of Proposition~\ref{geom}.
\end{proposition}
\begin{proof} Fix $x\in\R^n$ and $r>0$. If $A(x,r,2r)\subset B(0,1)$, then~\eqref{geom1} holds.
Suppose that there is a point $y\in A(x,r,2r)\setminus B(0,1)$. Let $y'=y+4ry/\abs{y}$. Note that
$B(y',2r)\cap B(0,1)=\varnothing$ and $2r\le \abs{y'-x}\le 6r$. For every line $L$ through $x$
there is a point $z\in B(y',r)$ such that $\dist(z,L)\ge r$. It follows that
$B(z,r/2)\cap S_{\alpha}(x,L)=\varnothing$ for sufficiently small absolute constant $\alpha$.
Therefore,
\[B(z,r/2)\subset A(x,r/4,8r)\setminus S_{\alpha}(x,L).\]
Furthermore,
$\mu(B(z,r/2))=\nu(B(z,r/2))$ because $B(z,r/2)\subset B(y',2r)$. Using the doubling property of $\nu$
we obtain
\[\nu(B(z,r/2))\ge C'\nu(B(x,2r))\ge \mu(A(x,r,2r))\]
where $C'$ depends on the doubling constant of $\nu$.
\end{proof}

\begin{example} The measure $\mu$ defined by the weight
\[w(x)=\begin{cases} \abs{x}^p,\quad \abs{x}> 1,\\
0,\quad \abs{x}\le 1,
\end{cases}\]
satisfies the assumptions of Theorem~\ref{dbthm2} whenever $-n<p<1-n$.
However, $\mu$ is not a doubling measure.
\end{example}

\begin{remark} For any $\gamma\in (n-1,n)$ and every Radon measure $\mu$
the Riesz potential $\Ri_{\gamma}\mu$ is comparable to $\norm{Df}$
for some quasiconformal mapping $f\colon\R^n\to\R^n$, provided that
$\Ri_{\gamma}\mu(x)\not\equiv\infty$~\cite[Lemma~4.1]{KM}.
By the composition formula for Riesz potentials~\cite[p.~118]{St70}
we have $\Ri_{\gamma}\mu=C \Ri_{n-1}\Ri_{\gamma-n+1}\mu$ for some constant $C=C(n,\gamma)$.
Since $\Ri_{\gamma-n+1}\mu$ is a doubling measure, and even an $A_1$-weight,
Lemma~4.1 in~\cite{KM} also follows from Theorem~\ref{dbthm}.
\end{remark}

\begin{remark} In~\cite{BHS06} Bonk, Heinonen, and Saksman used quasiconformal flows
to give another class of weights that are comparable to quasiconformal Jacobians.
These are weights of the form
\begin{equation}\label{bhsw}
w(x)=\exp\left\{-\int_{\R^n}\log\abs{x-z}\,d\mu(z)\right\},
\end{equation}
where $\mu$ is a signed Radon measure of sufficiently small total variation.
Neither of the classes of weights in Theorem~\ref{dbthm} and equation~\eqref{bhsw}
is contained in the other one.
\end{remark}

It is possible that the relation between Riesz potentials and Jacobians of quasiconformal mappings can
be extended beyond the results presented here.

\begin{question} Let $n\ge 2$, $0<\gamma<n$, and let $\mu$ be a Radon measure on $\R^n$ such that
$\Ri_{\gamma}\mu\not\equiv \infty$. Does there exist a quasiconformal mapping $f\colon\R^n\to\R^n$
such that $C^{-1}\Ri_{\gamma}\mu\le \det Df\le C\Ri_{\gamma}\mu$ a.e. for some constant $C$?
\end{question}

By Gehring's theorem~\cite[Theorem~1]{Ge73} every quasiconformal
mapping $f$ on $\R^n$ satisfies $\norm{Df}\in L_{\rm loc}^p$ for some $p>n$. This and
Theorem~\ref{dbthm} immediately imply the following result.

\begin{corollary}\label{RPI} If $\mu$ is a doubling measure on $\R^n$, $n\ge 2$, then the Riesz
potential $\Ri_{n-1}\mu$ is either infinite at every point or is in $L_{\rm loc}^p(\R^n)$ for some $p>n$.
\end{corollary}
Corollary~\ref{RPI} is probably known, although we could not find a reference.

\section{Isotropic doubling measures}\label{idm}

In this section we prove Theorem~\ref{dmidc} and give several examples of isotropic doubling measures.
It is clear that any weight $w\in L^{\infty}(\R^n)$ such that $\essinf_{x\in\R^n}w(x)>0$ is isotropic doubling.

\begin{lemma}\label{idc1} Suppose that a measure $\mu$ satisfies~\eqref{idmdef1}, and $n\ge 2$. Then:
\begin{enumerate}[(i)]
\item\label{idc1a} For any congruent rectangular boxes $R_1,R_2\subset \R^n$
\begin{equation}\label{idc11}
A^{-m}\le\frac{\mu(R_1)}{\mu(R_2)}\le A^m,\quad \text{where }
m=\left\lceil\frac{\dist(R_1,R_2)}{\diam(R_1)}\right\rceil+1.
\end{equation}
\item\label{idc1b} Let $Q\subset \R^n$ be a cube, and let $F$ be a face of $R$. The pushforward
$\pi_{\#}\mu_{|Q}$ of $\mu_{|Q}$
under the orthogonal projection $\pi\colon Q\to F$ is comparable to $\LL^{n-1}_{|F}$ with constants
that depend only on $n$ and $A$.
\end{enumerate}
\end{lemma}
\begin{proof} \eqref{idc1a} There exist congruent boxes $S_0,\dots, S_m$ such that $S_0=R_1$, $S_m=R_2$, and
$S_i\cap S_{i-1}\ne\varnothing$ for $i=1,\dots,m$. Repeated application of~\eqref{idmdef1} yields~\eqref{idc11}.

\eqref{idc1b} Let $l$ be the length of the edges of $Q$.
Let $Q_1$ and $Q_2$ be congruent $(n-1)$-dimensional cubes contained in $F$.
The boxes $R_i=\pi^{-1}(Q_i)$ have diameter at least $l$ and the distance between them
is at most $\diam F\le \sqrt{n-1}l$. By~\eqref{idc11} we have
\[
A^{-m}\le\frac{\pi_{\#}\mu(Q_1)}{\pi_{\#}\mu(Q_2)}\le A^m,\quad \text{where }
m=\left\lceil\sqrt{n-1}\right\rceil+1.
\]
The statement follows by letting $\diam(Q_i)\to 0$.
\end{proof}

Lemma~\ref{idc1} implies that any isotropic doubling measure on $\R^n$ vanishes on sets
of $(n-1)$-dimensional measure zero. The following lemma is the main step in our proof of Theorem~\ref{dmidc}.

\begin{lemma}\label{lemidc} Let $f : \R^n\to\R^n$ be a nonconstant $\delta$-monotone mapping, $n\ge 2$.
There exists $C>1$ such that for any rectangular box $R\subset\R^n$
\begin{equation}\label{lemidc1}
C^{-1}\frac{\diam (f(R))}{\diam (R)}\le
\frac{1}{\LL^n(R)}\int_R \norm{Df(x)}d\LL^n(x)\le C\frac{\diam (f(R))}{\diam (R)}.
\end{equation}
Here $C$ depends only on $\delta$ and $n$.
\end{lemma}
\begin{proof} Without loss of generality we may assume that
\[
R=\{x : \forall i \ 0\le x_i\le a_i\},\quad \text{where }\ a_1\ge a_2\ge \dots\ge a_n.
\]
We use $e_1,\dots,e_n$ to denote the standard basis vectors in $\R^n$.
Let $\pi(x_1,x_2,\dots,x_n)=(0,x_2,\dots,x_n)$ be the orthogonal projection onto the
linear span of $e_2,\dots,e_n$. Fix $z\in\pi(R)$ for now.
Since $f$ is quasisymmetric, it follows that
\[\abs{f(z+a_1e_1)-f(z)} \ge c\diam f(R)\]
where $c$ depends only on $\delta$ and $n$.
On the other hand, Definition~\ref{defdm} implies that the image of the line segment
$\pi^{-1}(z)\cap R$ under $f$ is a rectifiable curve of length at most
$\delta^{-1}\abs{f(z+a_1e_1)-f(z)}$, which is in turn bounded by $\delta^{-1}\diam f(R)$.
Therefore,
\begin{equation}\label{lemidc2}
c\diam f(R)\le \int_0^{a_1}\norm{Df(z+te_1)}\,dt \le \delta^{-1}\diam f(R).
\end{equation}
Integrating~\eqref{lemidc2} over $z\in\pi(R)$, we obtain
\[c\diam f(R)\prod_{i=2}^na_i\le
\int_R\norm{Df(x)}\,d\LL^n(x) \le \delta^{-1}\diam f(R)\prod_{i=2}^na_i.\]
Dividing the latter inequality by $\LL^n(R)$, we arrive at~\eqref{lemidc1}.
\end{proof}

\begin{proof}[Proof of Theorem~\ref{dmidc} ] By virtue of Lemma~\ref{lemidc} it suffices to
prove that the ratio $\diam f(R_1)/\diam f(R_2)$ is uniformly bounded for any congruent
boxes $R_1$ and $R_2$ with $R_1\cap R_2\ne\varnothing$.
Let $d=\diam R_1$ and pick a point $y\in R_1\cap R_2$.
There exists $z\in R_2$ such that $\abs{y-z}\ge d/2$. For any $x\in R_1$ we have
\[\frac{\abs{f(x)-f(y)}}{\abs{f(z)-f(y)}}\le \eta\left(\frac{\abs{x-y}}{\abs{z-y}}\right)\le\eta(2),\]
where $\eta$ is the modulus of quasisymmetry of $f$. Therefore,
\[\diam f(R_1)\le 2\eta(2)\abs{f(z)-f(y)}\le 2\eta(2)\diam f(R_2)\]
as required.
\end{proof}

Theorem~\ref{dmidc} provides a large family of isotropic doubling weights.
Indeed, for any set $E\subset \R^n$ of Hausdorff dimension $\dim E<1$ one can find
a $\delta$-monotone mapping $f\colon\R^n\to\R^n$ such that $\norm{Df}$
has essential limit $0$ at every point of $E$ (see Theorem~5.6~\cite{KM} and Theorem~17~\cite{Ko}).

Let $\mathcal{Q}$ and $\mathcal{R}$ denote the sets of arbitrarily oriented cubes
(resp. rectangular boxes) of diameter at most $1$. To a Radon measure $\mu$ on $\R^n$
we associate two maximal functions
\[\begin{split}
\mathcal{M}^{\mathcal{Q}}\mu(x)&=\sup\{\mu(Q)/\LL^n(Q)\colon x\in Q\in\mathcal{Q}\};\\
\mathcal{M}^{\mathcal{R}}\mu(x)&=\sup\{\mu(R)/\LL^n(R)\colon x\in R\in\mathcal{R}\}.
\end{split}\]
Obviously $\mathcal{M}^{\mathcal{Q}}\mu(x)\le \mathcal{M}^{\mathcal{R}}\mu(x)$ for all $x\in\R^n$.

\begin{lemma}\label{max} If $\mu$ is an isotropic doubling measure on $\R^n$,
then there exists $C\ge 1$ such that
$\mathcal{M}^{\mathcal{R}}\mu(x)\le C\mathcal{M}^{\mathcal{Q}}\mu(x)$ for all $x\in\R^n$.
\end{lemma}
\begin{proof} Fix $x\in\R^n$. Given a box $R\in \mathcal{R}$ containing $x$, find a cube $Q$
such that $R\subset Q$, the edges of $R$ and $Q$ are parallel, and the edges of $Q$ have the same
length as a longest edge of $R$. Although in general $Q\notin\mathcal{Q}$, the diameter of $Q$
cannot exceed $\sqrt{n}$. Lemma~\ref{idc1}~\eqref{idc1b} implies
$\mu(R)/\LL^n(R)\le C\mu(Q)/\LL^n(Q)$. The doubling property of $\mu$ and the definition of
$\mathcal{M}^{\mathcal{Q}}$ yield $\mu(Q)/\LL^n(Q)\le C\mathcal{M}^{\mathcal{Q}}(x)$.
\end{proof}

Lemma~\ref{max} and ~\cite[p.~5]{St70} imply that the maximal operator $\mathcal{M}^{\mathcal{R}}$
satisfies weak $(1,1)$ and strong $(p,p)$ estimates when restricted to the linear span of
isotropic doubling weights with the same doubling constant. In~\cite[Theorem~1.4]{AFN02} the operator
$\mathcal{M}^{\mathcal{R}}$ was proved to be bounded on certain Besov spaces. Neither Lemma~\ref{max}
nor Theorem~1.4~\cite{AFN02} imply each other.

\begin{example}\label{absx} The weight $w(x)=\abs{x}^p$, $x\in\R^n$, is isotropic doubling
if and only if $p>-1$.
Since the isotropic doubling condition is preserved under addition and weak${}^*$ convergence of measures,
for any Radon measure $\mu$ on $\R^n$ and any $\gamma\in (n-1,n)$ the
Riesz potential $\Ri_{\gamma}\mu$ is isotropic doubling unless $\Ri_{\gamma}\mu\equiv\infty$.
\end{example}
\begin{proof}
Let $\mu_p$ be the measure $\abs{x}^pd\LL^n(x)$, and let
\[\begin{split}
R_1&=\{x\in\R^n : \abs{x_1}\le 1; \ \forall i>1 \ \abs{x_i}\le \varepsilon\},\\
R_2&=\{x\in\R^n : \abs{x_1-2}\le 1; \ \forall i>1 \ \abs{x_i-1}\le \varepsilon\}.
\end{split}\]
If $-n<p\le -1$, then $\mu_p(R_1)/\varepsilon^{n-1}\to\infty$ as $\varepsilon\to0$.
Indeed, $\abs{x}\le \sqrt{n}x_1$ whenever $x_1\ge \abs{x_i}$ for all $i>1$.
The latter holds in particular when $x\in R_1$ and $x_1\ge \varepsilon$.
Therefore,
\[
\frac{\mu_p(R_1)}{\varepsilon^{n-1}}\ge n^{p/2}\int_{\varepsilon}^1 x_1^p dx_1\to\infty,
\quad \varepsilon\to0.
\]
Since $\mu_p(R_2)\le 2^n\varepsilon^{n-1}$, by Lemma~\ref{idc1}~\eqref{idc1b}
the measure $\mu_p$ is not isotropic doubling in the range $-n<p\le -1$.
The sufficiency part follows from Theorem~\ref{dmidc}.
Indeed, for $p>-1$ the mapping $f_p(x)=\abs{x}^{p}x$ is $\delta$-monotone by Proposition~16 in~\cite{Ko}.
By Theorem~\ref{dmidc} $\norm{Df}$ is an isotropic doubling weight.
Since $\norm{Df_p}$ is a $p$-homogeneous radially symmetric function, it coincides with
a constant multiple of $\abs{x}^p$.
\end{proof}

Given two points $a,b\in\R^n$,
let $[a,b]$ denote the segment $\{\lambda a+(1-\lambda)b\colon \lambda\in [0,1]\}$.
A set $A\subset\R^n$ is \emph{uniformly linearly non-convex} (ULNC)
if there exists $\tau>0$ such that the following holds:
To each pair of distinct points $a,b\in A$ there corresponds $c\in [a,b]$ such that
$B(c,\tau\abs{a-b})\cap A=\varnothing$.~\cite{TW05} The following lemma contains
an equivalent definition of ULNC sets.

\begin{lemma}\label{ULNCb} Suppose that $A\subset\R^n$ is ULNC with a constant $\tau>0$.
Let $\tau'=\tau/(2+2\tau)$. Then to each pair of distinct points $a,b\in \R^n$ there
corresponds $c\in [a,b]$ such that $B(c,\tau'\abs{a-b})\cap A=\varnothing$.
\end{lemma}
\begin{proof} Without loss of generality we may assume that $\abs{a-b}=1$.
Suppose that there are points $a',b'\in A$ such that
$\abs{a'-a}\le \tau'$ and $\abs{b'-b}\le \tau'$; otherwise we can set $c=a$ or $c=b$.
Since $\abs{a'-b'}\ge (1-2\tau')$, there exists
$c'\in [a',b']$ such that $B(c',\tau(1-2\tau'))$ is disjoint from $A$.
Each point of the segment $[a',b']$ is at distance at most $\tau'$ from a point of $[a,b]$.
Therefore, there exists $c\in [a,b]$ such that $\abs{c-c'}\le\tau'$. The closed ball
with center $c$ and radius
\[\tau(1-2\tau')-\tau'=\tau\left(1-\frac{\tau}{1+\tau}\right)-\frac{\tau}{2+2\tau}=
\frac{\tau}{2+2\tau}=\tau'\]
is disjoint from $A$, as required.
\end{proof}

\begin{proposition}\label{ULNC} Let $A\subset\R^n$ be a nonempty ULNC set. Given $p>0$, define
$w(x)=\dist(x,A)^p$ for $x\in\R^n$. The weight $w$ is isotropic doubling.
\end{proposition}
The proof is preceded by a lemma.

\begin{lemma}\label{cts} A continuous weight $w\in C(\R^n)$ is isotropic doubling if and only
if there exists a constant $C$ such that
\begin{equation}\label{cts1}
\int_0^1 w(x_0+tv)\,dt \le C\int_0^1 w(x_0+tv')\,dt
\end{equation}
for all points $x_0\in \R^n$ and all vectors $v,v'\in\R^n$ such that $\abs{v}=\abs{v'}>0$.
Moreover, the constant $A$ in~\eqref{idmdef1} depends only on $C$ in~\eqref{cts1}.
\end{lemma}
\begin{proof} If $w$ is continuous and isotropic doubling, then~\eqref{cts1} follows by
applying~\eqref{idmdef1} to two thin rectangular boxes containing the segments
$\{x_0+tv\colon 0\le t\le 1\}$ and $\{x_0+tv'\colon 0\le t\le 1\}$.
Conversely, suppose that~\eqref{cts1} holds. Let $R_1$ and $R_2$ be congruent boxes with nonempty
intersection. Let $d$ be the length of a longest edge of $R_i$. If $L_i\subset R_i$ is a segment
of length $d$ connecting two opposite faces of $R_i$, then $R_i$ is contained in the closed
$\sqrt{n-1}d$-neighborhood of $L_i$ for each $i\in{1,2}$. This implies $\dist(L_1,L_2)\le 2\sqrt{n-1}d$.
For an integer $N\ge 2\sqrt{n-1}+2$ we can find a sequence of segments $l_1,\dots,l_N$ of equal length such that
$l_1=L_1$, $l_N=L_2$, and $l_i$ shares an endpoint with $l_{i+1}$ for $i=1,\dots,N-1$. A repeated
application of~\eqref{cts1} yields
\[\int_{L_1} w(x)\,ds \le C^{N-1}\int_{L_2} w(x)\,ds.\]
Since the box $R_i$ is foliated by segments such as $L_i$, inequality~\eqref{idmdef1} holds
with $A=C^{N-1}$.
\end{proof}

\begin{proof}[Proof of Proposition~\ref{ULNC} ] If $\dist(x_0,A)\ge 2\abs{v}$, then
\[
\dist(x_0,A)\le \left(\int_0^1 w(x_0+tv)\,dt\right)^{1/p}\le 3\dist(x_0,A).
\]
If $\dist(x_0,A)\le 2\abs{v}$, then
\[
\left(\int_0^1 w(x_0+tv)\,dt\right)^{1/p}\le 3\abs{v}.
\]
Using Lemma~\ref{ULNCb} we obtain the reverse inequality
\[
\left(\int_0^1 w(x_0+tv)\,dt\right)^{1/p}\ge c\abs{v},
\]
where $c>0$ depends only on $\tau$ in the definition of a ULNC set.
By Lemma~\ref{cts} the weight $w$ is isotropic doubling.
\end{proof}

\begin{remark} Theorem~\ref{dmidc} admits no converse, that is, there exist isotropically doubling weights
that are not comparable to $\norm{Df}$ for any $\delta$-monotone mapping $f$.
This follows from the results of Semmes~\cite{Se96} and Laakso~\cite{La02}, who constructed compact ULNC sets $A\subset \R^n$
such that for some $p>0$ the weight $w(x):=\dist(x,A)^p$ is not comparable to $\norm{Df}$
for any quasiconformal mapping $f\colon\R^n\to\R^n$. By Proposition~\ref{ULNC} the weight $w$ is isotropic doubling.
The example of Semmes is a form of the Antoine's
necklace, and requires $n\ge 3$ (see the proof of Theorem~1.10 in~\cite{Se96}). The example
of Laakso is two-dimensional~\cite[Theorem~1.7]{La02}. The fact that the sets considered
by Semmes and Laakso are ULNC follows from Theorem~3.3 in~\cite{TW05}.
\end{remark}

\begin{remark} Buckley, Hanson and MacManus considered several stronger versions of the
doubling condition in~\cite{BHM01}, but our Definition~\ref{idmdef}
does not appear to be directly related to any of them.
\end{remark}

\section{The Monge-Amp\`ere measures}\label{quconvex}

The main purpose of this section is to point out a consequence of
Theorem~\ref{dmidc} for convex functions on $\R^n$.
Let $u\colon\R^n \to \R$ be a convex function, $n\ge 1$.
The \emph{subdifferential} of $u$ at a point $z\in \R^n$ is the set
\[
\partial u(z) = \{ p \in \R^n\colon u(x) \ge u(z) + \inn{p, x-z} \ \forall x \in \R^n\}.
\]
The convexity of $u$ implies $\partial u(z)\ne \varnothing$ for all $x\in\R^n$.
The function $u$ is differentiable at $z\in\R^n$ if and only if $\partial u(z)$ is the
one-point set $\{\nabla u(z)\}$. Given $z \in \R^n$ and $p \in \partial u(z)$, let
\[
u_{z,p}(x)= u(x) -u(z) - \inn{p, x-z}, \quad x\in\R^n.
\]
The \emph{section}~\cite[p.~45]{Gu01} of $u$ with the center $z\in\R^n$,
direction $p \in \partial u(z)$, and height $t>0$ is defined as
\[
S_u(z,p,t)=\{x \in \R^n\colon u_{z,p}(x)< t\}.
\]
If $u$ is G\^ateaux differentiable at $z$, then we write $u_z$ instead of $u_{z,\nabla u(z)}$.
The \emph{Monge-Amp\`ere measure} associated with $u$ is defined
by $\mu_u(E)=\mathcal{L}^n(\partial u(E))$ for every Borel set
$E\subset\R^n$. Note that adding an affine function to $u$ does not change $\mu_u$.
If $u\in C^2$, then $\mu_u$ is absolutely continuous with the density $\det D^2 u(x)$,
where $D^2u$ is the Hessian matrix of $u$.

\begin{definition}\label{rounddef}~\cite{KM}
A convex function $u\colon\R^n\to\R$ has \emph{round sections} if there exists
a constant $\tau\in(0,1)$ such that for every $z \in \R^n$, $p \in \partial u(z)$ and $t>0$
\begin{equation}\label{round}
B(z,\tau R)\subset S_u(z,p,t)\subset B(z,R)
\end{equation}
for some $R\in (0,\infty)$.
\end{definition}

\begin{proposition} If a convex function $u\colon \R^n\to\R$ has round sections, then
$||D^2u||$ is an isotropic doubling weight.
\end{proposition}
\begin{proof} First, consider the case $n=1$. By Theorem~3.1~\cite{KM} and
Remark~3.3~\cite{KM} $u$ is continuously differentiable and $\nabla u\colon\R\to\R$
is quasisymmetric. If $I$ and $J$ are adjacent intervals of equal length, then
$\mu_u(I)\le \eta(1)\mu_u(J)$, where $\eta$ is as in~\ref{def3a}.
Thus $\mu_u$ is doubling, which in one dimension is the same as being isotropic
doubling.

Next, assume $n\ge 2$. By Theorem~3.1~\cite{KM} $u$ is continuously differentiable, and the gradient
mapping $\nabla u\colon\R^n\to\R^n$ is quasiconformal. Since the change of variables
formula applies to quasiconformal mappings~\cite[Theorem~33.3]{Va71}, we have
\[\mu_u(E)=\int_E \det D^2u(x)\,d\LL^n(x),\quad E\subset\R^n \text{ measurable.}\]
Quasiconformality also implies that $\det D^2u$ is comparable to $\norm{D^2 u}^n$ up to a multiplicative
constant. By Theorem~14~\cite{Ko} the gradient $\nabla u$ is a $\delta$-monotone mapping,
hence $\norm{D^2 u}$ is isotropic doubling by Theorem~\ref{dmidc}.
\end{proof}

\begin{remark}\label{pointless} The mapping $f_{\mu}$ defined by~\eqref{db0a} is the gradient of the convex function
\begin{equation}\label{vmu}
v_{\mu}(x):=\int_{\R^n}\left\{\abs{x-z}-\abs{z}+\frac{\inn{x,z}}{\abs{z}}\right\}\,d\mu(z).
\end{equation}
Indeed, for all $x,z\in\R^n$
\[\frac{x-z}{\abs{x-z}}+\frac{z}{\abs{z}}\in \partial\left(\abs{x-z}-\abs{z}+\frac{\inn{x,z}}{\abs{z}}\right),\]
where the subgradient is taken with respect to $x$. Therefore, $f_{\mu}(x)\in\partial v_{\mu}(x)$
for all $x\in\R^n$. The continuity of $f_{\mu}$ implies that $v_{\mu}(x)$ is differentiable and
$\nabla v_{\mu}=f_{\mu}$.
\end{remark}

If the gradient of a convex function is quasisymmetric, then
it is $\delta$-monotone by Lemma~13 in~\cite{Ko}. This together with
Remarks~\ref{close} and~\ref{pointless} indicate that condition~\eqref{dbthm3} is close to being best
possible for the quasiconformality of $f_{\mu}$ in Theorem~\ref{dbthm2}.

\section{Singular isotropic doubling measures}\label{singular}

It is well-known that doubling measures on $\R^n$ can be supported on sets of arbitrarily small
positive Hausdorff measure (\cite{BA56} and~\cite[p.~40]{St93}). This is not true for
isotropic doubling measures in dimensions $n\ge 2$ by Lemma~\ref{idc1}~\eqref{idc1b}.
In this section we prove Theorem~\ref{sing}, which
 shows that isotropic doubling measures can be singular with respect
to $\LL^n$. We will use the following notation: $\# A$ is the cardinality of a finite set $A$,
$a\vee b=\max\{a,b\}$, $a\wedge b=\min\{a,b\}$.

\begin{proof}[Proof of Theorem~\ref{sing} ] Let us introduce a sequence of vectors
\[w^k:=4^{nk^2}(4^k,4^{2k},4^{3k},\dots,4^{nk})\in\R^n, \quad k=1,2,\dots.\]
For $k\ge 1$ define
\[\lambda_k(x)=1+\frac12\cos(\inn{x,w^k}),\quad x\in\R^n.\]
and
\[\Lambda_m(x)=\prod_{k=1}^m \lambda_k(x).\]
We shall obtain a singular isotropic doubling measure $\mu$ as the weak${}^*$ limit of a subsequence
of $\{\Lambda_m(x)d\LL^n(x)\colon m\ge 1\}$. The proof involves several steps.

\emph{Step $1$. Existence of $\mu$.} For each $m$ the weight $\Lambda_m$ is $2\pi$-periodic in all
variables. Let $Q=[-\pi,\pi]^n$. Clearly $\int_Q \Lambda_1(x)\,d\LL^n(x)=(2\pi)^n$.
For $m\ge 2$ we obtain $\int_Q (\Lambda_m-\Lambda_{m-1})d\LL^n=0$ by integrating
in $x_n$ over the interval $[-\pi,\pi]$ and using the orthogonality of the trigonometric family. Therefore,
$\int_Q \Lambda_m\,d\LL^n=(2\pi)^n$ for all $m\ge 1$. Choose a subsequence of $\Lambda_m(x)\,d\LL^n(x)$
converging in the weak${}^*$ sense to a $(2\pi)$-periodic measure $\mu$ with $\mu(Q)=(2\pi)^n$.

\emph{Step $2$. Lipschitz estimates.}
Since the Euclidean norm of $w^k$ satisfies
\begin{equation}\label{normwk}
4^{n(k^2+k)}\le \abs{w^k}\le 2\cdot 4^{n(k^2+k)},
\end{equation}
$\lambda_k$ is a Lipschitz function with the constant $L_k:=4^{n(k^2+k)}$.
Therefore, for any $x,x'\in\R^n$ we have
\begin{equation}\label{Lip1}
\frac{\lambda_k(x')}{\lambda_k(x)}\ge \frac{1/2}{1/2+L_k\abs{x-x'}}\ge 1-2L_k\abs{x-x'}.
\end{equation}
Multiplying~\eqref{Lip1} over $k=1,\dots,m$ yields
\begin{equation}\label{Lip2}
\frac{\Lambda_m(x')}{\Lambda_m(x)}\ge \prod_{k=1}^m(1-2L_k\abs{x-x'})\ge 1-2\abs{x-x'}\sum_{k=1}^m L_k
\ge 1-3L_m\abs{x-x'}
\end{equation}
provided that $3L_m\abs{x-x'}<1$. Thus
\begin{equation}\label{Lip3}
1-3L_m\abs{x-x'}\le \frac{\Lambda_m(x')}{\Lambda_m(x)}\le (1-3L_m\abs{x-x'})^{-1},\quad \abs{x-x'}<\frac{1}{3L_m}.
\end{equation}

\emph{Step $3$. Singularity of $\mu$.} Given $x\in\R^n$ and $r>0$, let $Q(x,r)$ denote the cube
$\{y\in\R^n\colon \forall i\ \abs{x_i-y_i}\le r\}$.
To prove that $\mu$ is purely singular, it suffices to show that
\begin{equation}\label{to0}
\liminf_{r\to0}r^{-n}\mu(Q(x,r))=0\quad\text{for a.e. }x\in\R^n.
\end{equation}
For a.e. $x\in \R^n$ we have $\lim_{m\to\infty}\Lambda_m(x)=0$ by~\cite[p.~209]{Zy02}. Fix such a point $x$.
For $m\ge 1$ let $r_m=4^{-n(m+1)^2}\pi$. Note that the functions $\lambda_k$, $k>m$, are $2r_m$-periodic in each variable.
Using orthogonality as in Step~$1$, we obtain
\[\mu(Q(x,r_m))=\int_{Q(x,r_m)}\Lambda_m(x)\,d\LL^n(x).\]
For every $y\in Q(x,r_m)$ we have
\[3L_m\abs{x-y}\le 3\cdot 4^{n(m^2+m)}\sqrt{n}4^{-n(m+1)^2}\pi\to 0\quad\text{as }m\to\infty.\]
Therefore, $\sup\{\Lambda_m(y) \colon y\in Q(x,r_m)\}\to 0$ as $m\to\infty$. This proves~\eqref{to0}.

It remains to prove that $\Lambda_m$ satisfies~\eqref{idmdef1} with a constant independent of $m$.
Consider a segment $L$ given by the parametric equations $x=x_0+tv$, $0\le t\le 1$,
of length $r=\norm{v}>0$. Set $K=\max\{k\ge 1\colon 4^{nk^2} \le 1/r\}$ or $K=0$ if no such $k$ exists.
Our goal is to show that
\begin{equation}\label{ourgoal}
\int_{L} \Lambda_m\,ds \underset{n}\approx r\prod_{k=1}^{K\wedge m}\left(1+\frac12\cos(\inn{x_0,w^k})\right)
\end{equation}
where $\underset{n}{\approx}$ indicates that multiplicative constants in the estimate depend
only on $n$. The estimate~\eqref{ourgoal} implies~\eqref{idmdef1} by Lemma~\ref{cts}.

\emph{Step $4$. Low-frequency terms.}
If $K\ge 2$, then the terms $\lambda_k$ with $1\le k\le K-1$ are essentially constant
on $L$. Indeed,~\eqref{normwk} implies
\begin{equation}\label{small}
\abs{\inn{v,w^{k}}}\le 2r4^{n(k^2+k)}\le 2\cdot 4^{-nK}.
\end{equation}
This implies that for all $0\le t\le 1$ and all $1\le k\le K-1$ we have
\[\abs{\lambda_k(x_0+tv)-\lambda_k(x_0)}\le 4^{-nK},\]
hence
\[1-2\cdot 4^{-nK}\le \frac{\lambda_k(x_0+tv)}{\lambda_k(x_0)}\le 1+2\cdot 4^{-nK}.\]
Taking a product over $k=1,\dots,K-1$ yields
\begin{equation}\label{low1}
(1-2\cdot 4^{-nK})^{K-1}\le
\prod_{k=1}^{K-1}\frac{\lambda_k(x_0+tv)}{\lambda_k(x_0)} \le (1+2\cdot 4^{-nK})^{K-1}.
\end{equation}
Since
\[
(1-2\cdot 4^{-nK})^{K-1}\ge 1-2(K-1)4^{-nK}\ge \frac{1}{2}
\]
and
\[
(1+2\cdot 4^{-nK})^{K-1}\le (1-2\cdot 4^{-nK})^{1-K}\le 2,
\]
inequality~\eqref{low1} implies
\begin{equation}\label{mainint2}\begin{split}
\frac{1}{2}\int_{L} \Lambda_m\,ds &\le r \prod_{k=1}^{m\wedge (K-1)}\left(1+\frac12\cos(\inn{x_0,w^k})\right)
\\
&\times\int_0^1 \prod_{k=m\wedge (K-1)}^m \left(1+\frac12\cos(\inn{x_0+tv,w^k})\right)\,dt \le 2\int_{L} \Lambda_m\,ds.
\end{split}\end{equation}
The estimate~\eqref{mainint2} proves~\eqref{ourgoal} in the case $m\le K+1$.

\emph{Step $5$. High-frequency terms.} For the rest of the proof we assume $m\ge K+2$.
Let $v^*=(\abs{v_1},\abs{v_2},\dots,\abs{v_n})$ be the vector whose components are the absolute values
of the components of $v$. Introduce two sets of integers depending on $v$:
\[\begin{split}
G&=\left\{K+2 \le k\le m \colon \abs{\inn{v,w^k}}\ge\frac14\inn{v^*,w^k}\right\};\\
B&=\left\{K+2 \le k\le m \colon k\notin G\right\}.
\end{split}\]
For every $k\in G$ the restriction of $\lambda_k$ to $L$ is rapidly oscillating, because
\begin{equation}\label{large}
\abs{\inn{v,w^k}}\ge 4^{nk^2-1}\sum_{i=1}^n4^{ik}\abs{v_i}
\ge 4^{nk^2+k-1}\abs{v}\ge 4^{nk^2+k-1}4^{-n(K+1)^2}\ge 4^{nk}.
\end{equation}
Furthermore, if $k\in G$ and $k>l\ge 1$, then the restriction of $\lambda_k$ to $L$
has period at least $4^{nk}$ times smaller than the period of the restriction of $\lambda_l$:
\begin{equation}\label{riesz1}
\abs{\inn{v,w^k}}\ge 4^{nk^2-1}\sum_{i=1}^n 4^{ik}\abs{v_i}\ge 4^{n(k^2-l^2)-1}
\sum_{i=1}^n4^{nl^2+il}\abs{v_i}\ge 4^{nk}\abs{\inn{v,w^l}}.
\end{equation}
We have chosen the vectors $w^k$ so that $\# B$ is bounded by a constant depending only on $n$,
by Lemma~\ref{domination} below. Hence
\begin{equation}\label{mainint3}
\prod_{k=K}^m \lambda_k(x)\underset{n}{\approx} \prod_{k\in G}\lambda_k(x),\quad x\in\R^n.
\end{equation}
Observe that
\begin{equation}\label{mainint9}
\prod_{k\in G}\left(1+\frac12\cos(\inn{x_0+tv,w^k})\right)
=1+\sum_{A\subseteq G, \ A\ne\varnothing}\frac{1}{2^{\# A}}\prod_{k\in A}\cos(\inn{x_0+tv,w^k}).
\end{equation}
Converting the trigonometric product into a sum, we obtain
\[\prod_{k\in A}\cos(\inn{x_0+tv,w^k})
=\frac{1}{2^{\# A}}\sum_{\pm} \cos\left(\sum_{k\in A} \pm \inn{x_0+tv,w^k}\right),
 \]
where the exterior sum is taken over all $2^{\# A}$ choices of the $\pm$ signs.
Using~\eqref{riesz1}, we find that for any choice of the signs
\[\Abs{\sum_{k\in A} \pm \inn{v,w^k}}\ge \frac{2}{3}\abs{\inn{v,w^{\max A}}},\]
where $\max A$ is the maximal element of $A$. Therefore,
\[\Abs{\int_0^1 \prod_{k\in A}\cos(\inn{x_0+tv,w^k})\,dt}
\le 3\abs{\inn{v,w^{\max A}}}^{-1} \le 3\cdot 4^{-n\max A},\]
where the last inequality follows from~\eqref{large}.
For any $l\in G$ there exist at most $2^{l-K-2}$ nonempty subsets $A\subseteq G$ such that $\max A=l$.
Thus we can estimate the righthand side of~\eqref{mainint9} as follows.
\begin{equation}\label{mainint10}\begin{split}
\sum_{A\subseteq G, \ A\ne\varnothing}\frac{1}{2^{\# A}}\Abs{\int_0^1 \prod_{k\in A}\cos(\inn{x_0+tv,w^k})\,dt}
&\le 3\sum_{l=K+2}^{m} 4^{-nl}2^{l-K-2}\\ &\le 3\sum_{l=K+2}^{\infty} 4^{-l}
=4^{-K-1}\le\frac{1}{4}.
\end{split}\end{equation}
Combining~\eqref{mainint9} and~\eqref{mainint10} yields
\[\frac{3}{4}\le \int_0^1 \prod_{k\in G}\left(1+\frac12\cos(\inn{x_0+tv,w^k})\right)\,dt\le\frac{5}{4},\]
which together with~\eqref{mainint2} and~\eqref{mainint3} imply
\begin{equation}\label{mainint0}
\int_{L} \Lambda_m\,ds \underset{n}{\approx} r \prod_{k=1}^{K}\left(1+\frac12\cos(\inn{x_0,w^k})\right).
\end{equation}
This proves~\eqref{ourgoal}, and thereby the theorem.
\end{proof}

\begin{remark} The above construction is not symmetric with respect to the variables $x_1,\dots,x_n$.
This feature of $\mu$ will be used to prove Theorem~\ref{nobl}. But one can construct
a more symmetric singular measure $\nu=\sum_{\sigma} T^{\sigma}_{\#}\mu$, where
$T^{\sigma}(x_1,\dots,x_n)=(x_{\sigma(1)},\dots,x_{\sigma(n)})$ and the summation is over all
permutations of $n$ elements. The isotropic doubling condition is obviously preserved under
addition of measures.
\end{remark}

Now we prove Lemma~\ref{domination}, which was used in Step~5 of the proof of Theorem~\ref{sing}
with parameters $q=4$ and $\varepsilon=3/4$.

\begin{lemma}\label{domination} Let $(v_1,\dots, v_n)$ be a nonzero vector in $\R^n$, where $n\ge 2$.
Given $q>1$ and $\varepsilon\in (0,1)$, define
\[B(q,\varepsilon)=\left\{k\in\Z \colon \Abs{\sum_{i=1}^n q^{ik}v_i}<\frac{1-\varepsilon}{1+\varepsilon}
\sum_{i=1}^n q^{ik}\abs{v_i}\right\}.\]
Then
\begin{equation}\label{dom1}
\# B(q,\varepsilon)\le n(n-1)\frac{\log((n-1)/\varepsilon)}{\log q}.
\end{equation}
\end{lemma}
\begin{proof} We claim that
\begin{equation}\label{dom11}
B(q,\varepsilon) \subseteq \bigcup_{1\le i<j\le n} A_{ij},
\end{equation}
where
\[
A_{ij}=\left\{k\in\Z\colon \frac{\varepsilon}{n-1} \le \frac{q^{ik}\abs{v_i}}{q^{jk}\abs{v_j}}\le \frac{n-1}{\varepsilon}\right\}
\]
if $v_j\ne 0$ and $A_{ij}=\varnothing$ otherwise.
Indeed, suppose that $k\in\Z\setminus\bigcup_{1\le i<j\le n} A_{ij}$. Let $m$ be such that
$q^{ik}\abs{v_i}$ is maximal when $i=m$. Then
\[
\Abs{\sum_{i=1}^n q^{ik}v_i}\ge q^{mk}\abs{v_m}-\sum_{i\ne m}q^{ik}\abs{v_i}
\ge (1-\varepsilon)q^{mk}\abs{v_m},
\]
while on the other hand
\[
\sum_{i=1}^nq^{ik}\abs{v_i}\le q^{mk}\abs{v_m}+(n-1)\frac{\varepsilon}{n-1}q^{mk}\abs{v_m}
=(1+\varepsilon)q^{mk}\abs{v_m}.
\]
Thus $k\ne B(q,\varepsilon)$, which proves~\eqref{dom11}. It remains to estimate the cardinality
of the sets $A_{ij}$, $1\le i<j\le n$. We may assume that $v_j\ne 0$, for otherwise $A_{ij}$ is empty.
Since
\[
A_{ij}=\left\{k\in\Z\colon \log \frac{\varepsilon}{n-1} \le
k(i-j)\log q+\log\frac{\abs{v_i}}{\abs{v_j}}\le \log\frac{n-1}{\varepsilon}\right\},
\]
it follows that
\[
\# A_{ij} \le \frac{2\log((n-1)/\varepsilon)}{\abs{i-j}\log q}\le \frac{2\log((n-1)/\varepsilon)}{\log q},
\]
which after summation yields~\eqref{dom1}.
\end{proof}

\begin{question}\label{absq} Is it true that every isotropic doubling measure on $\R^n$ is
absolutely continuous with respect to the $s$-dimensional Hausdorff measure for all $s<n$?
\end{question}
Lemma~\ref{idc1}~\eqref{idc1b} gives an affirmative answer to Question~\ref{absq} in the case $s\le n-1$.

\section{Bi-Lipschitz transformations}\label{BL}

A mapping $f : \R^n\to\R^n$ is called bi-Lipschitz if there exists $L\ge 1$ such that
\[
L^{-1}\abs{x-y}\le \abs{f(x)-f(y)}\le L\abs{x-y}
\]
for all $x, y\in \R^n$. It is easy to see that the doubling condition is invariant under
bi-Lipschitz mappings. In contrast, Theorem~\ref{nobl} asserts that the isotropic doubling condition
is not bi-Lipschitz invariant in dimensions $n\ge 2$. Its proof is preceded by a definition and a lemma.

\begin{definition} An unbounded Jordan curve $\Gamma \subset \R^2$ is a \emph{chord-arc curve} if
there is a constant $C$ such that for any two points $a,b\in\Gamma$ the part of $\Gamma$
bounded by $a$ and $b$ has length at most $C\abs{a-b}$.
\end{definition}

Given a set $A\subset\R^n$ and $\varepsilon>0$, we define
$N_{\varepsilon}A=\{x\in\R^n\colon \dist(x,A)\le\varepsilon\}$.

\begin{lemma}\label{charc} Suppose that $\mu$ is a Radon measure on $\R^n$, $n\ge 2$ such that
for every bi-Lipschitz mapping $f\colon\R^n\to\R^n$ the pushforward measure
$f_{\#}\mu$ is isotropic doubling. Then for any chord-arc curve $\Gamma\colon\R\to\R^2\subset\R^n$
and any bounded subarc $\Gamma'\subset \Gamma$ we have
\begin{equation}\label{charc1}
0<\liminf_{\varepsilon\to 0} \varepsilon^{1-n}\mu(N_{\varepsilon}\Gamma')\le
 \limsup_{\varepsilon\to 0} \varepsilon^{1-n}\mu(N_{\varepsilon}\Gamma')<\infty.
\end{equation}
\end{lemma}
\begin{proof} There exists a bi-Lipschitz mapping $g\colon\R^2\to\R^2$ that transforms $\Gamma$
into the real axis, see~\cite[p.~89]{Tu81} or~\cite[Proposition~1.13]{JK82}. Let $f\colon\R^n\to\R^n$
be the bi-Lipschitz extension of $g$ that acts trivially on the remaining $(n-2)$ coordinates.
Let $L$ be a bi-Lipschitz constant of $f$. Note that
\[N_{\varepsilon/L}f(\Gamma') \subset f(N_{\varepsilon}\Gamma')\subset N_{L\varepsilon}f(\Gamma'),\]
and $f(\Gamma')$ is a line segment. Lemma~\ref{idc1}~\eqref{idc1b} implies that for any line segment $I$
the quantity $\varepsilon^{1-n}\mu(N_{\varepsilon}I)$ remains bounded away from $0$ and $\infty$
as $\varepsilon\to0$.
\end{proof}

\begin{proof}[Proof of Theorem~\ref{nobl} ] We use the notation of Theorem~\ref{sing}.
Define a sequence of nonnegative functions on $\R$ as follows: $h_0(t)=-t$ for all $t\in\R$,
\begin{equation}\label{nobl1}
h_k(t)=\min\{s\ge h_{k-1}(t)\colon \lambda_k(t,s,0,\dots,0)=3/2\}\quad\text{for }k\ge 1.
\end{equation}
Note that $\lambda_k(t,s,0,\dots,0)=3/2$ precisely when $\cos(4^{nk^2+k}(t+4^ks))=1$.
Therefore, $h_k$ is a decreasing piecewise affine function such that $h_k'(t)=-4^{-k}$ at the points where $h_k$
is continuous. Also,
\begin{equation}\label{nobl2}
0\le h_k(t)-h_{k-1}(t)\le 2\pi 4^{-nk^2-2k},\quad t\in\R.
\end{equation}
Therefore, the limit $h=\lim_{k\to\infty}h_k$ is a decreasing function on $\R$.
Using~\eqref{nobl2} and the Lipschitz property of $\lambda_k$, we find that for any $1\le k<l$
we have
\begin{equation}\label{nobl3}\begin{split}
\lambda_k(t,h_l(t),0,\dots,0)&\ge \lambda_k(t,h_k(t),0,\dots,0)-4^{nk^2+2k}(h_l(t)-h_k(t))\\
&\ge \frac{3}{2}-2\pi 4^{nk^2+2k}\sum_{j=k+1}^{l}4^{-nj^2-2j}
\ge \frac{3}{2}-4^{-2nk}.
\end{split}
\end{equation}
Therefore, $\Lambda_k(t,h(t),0,\dots,0)\to\infty$ as $k\to\infty$ uniformly in $t\in\R$.
Let $\Gamma$ be the curve obtained from the graph of $h$ by filling the gaps
at the points of discontinuity with vertical segments. Since $h$ is a decreasing function,
$\Gamma$ is a chord-arc curve. Let $G'=G\cap \{(x,y)\colon \abs{x}\le 1\}$.

In view of Lemma~\ref{charc} it remains to prove that
\begin{equation}\label{neig}
\limsup_{\varepsilon\to0}\varepsilon^{-1}\mu(N_{\varepsilon}\Gamma')=\infty
\end{equation}
Pick a large positive integer $K$ and let $\varepsilon=4^{-nK^2}$. For any line segment $L$ of
length $\varepsilon$ we have
\begin{equation}\label{larrge}
\int_L \Lambda_m\,ds\underset{n}\approx \varepsilon\max_L\Lambda_K,\quad m\ge K,
\end{equation}
by virtue of~\eqref{ourgoal}. For segments $L$ at distance less than $\varepsilon$ from the graph
of $h$, the righthand side of~\eqref{larrge} is large by the Lipschitz estimate~\eqref{Lip3}
applied to $\Lambda_{K-1}$. Letting $K\to\infty$ we obtain~\eqref{neig}.
\end{proof}

\section*{Acknowledgements}

We thank Robert Kaufman and Jeremy Tyson for helpful discussions.
We also thank the referee for several comments and corrections.

\bibliographystyle{amsalpha}

\end{document}